\magnification1200

\catcode`?=11 \relax%
\gdef\XA?X{?}%
\gdef\?empty{}%
\long\gdef\Thm#1#2{%
\noindent{\bf#1.\ }{\sl\ignorespaces#2}\medskip
}%
\long\gdef\Def#1#2{%
\noindent{\bf#1.\ }{\rm\ignorespaces#2}\medskip
}%
\global\let\Itemitem=\item%
\gdef\Matrix#1{\gdef\CR{\cr}\matrix{#1}}%
\gdef\pf{\par\noindent{\sl Proof:}\ }%
\gdef\QED{\penalty 10000 \hskip 3pt plus 1.6pt minus 1pt%
\vrule height4 true pt width3 true pt depth0pt %
\ifmmode\else\par\bigskip\fi%
}
\global\let\Qed=\QED
\gdef\Frac#1#2{{{#1}\over{#2}}}%
\gdef\NumAlign#1{%
\def\CR##1{\gdef\x??x{##1}%
\ifx\x??x\?empty\def\x??x{\cr}%
\else\def\x??x{&##1\cr}\fi\x??x}%
$$\leqalignno{#1}$$
}%
\global\font\SectionFont=cmbx12%
\newcount\sectionnumber\global\sectionnumber=1 %
\gdef\section#1{\bigskip\penalty-1000\noindent{{\SectionFont%
\number\sectionnumber. #1}\penalty10000\medskip}%
\global\advance\sectionnumber by 1 }%

\catcode`?=11 \relax%
\newbox\arrowtmp
\newbox\arrowtmpA
\newdimen\arrowwd

\def\limitsarrow#1#2#3{%
\setbox\arrowtmp=\hbox{$\scriptstyle11#1$}%
\setbox\arrowtmpA=\hbox{$\scriptstyle11#2$}%
\ifdim\wd\arrowtmp<\wd\arrowtmpA\relax
\arrowwd=\wd\arrowtmpA\relax\else%
\arrowwd=\wd\arrowtmp\relax\fi%
\smash{\mathop{\vbox to 3pt{\vss\hbox to \arrowwd{#3}}}\limits^{#1}\limits_{#2}}%
}

\def\RA#1{%
\limitsarrow{#1}{}{\rightarrowfill}}
\catcode`@=11

\def\mbig#1#2{{\hbox{$\left#1\vtop to#2{\vfil}\right.\n@space$}}}

\def\Downarrow#1#2{\arrowwd=#1 \advance\arrowwd by#2 %
\raise#1\vtop to0pt{\hsize=0pt\mbig{\downarrow}%
{\arrowwd}\vss}\vtop to#2{\hsize=0pt\vfil}%
\advance\arrowwd by #1 \advance\arrowwd by 6pt%
\vbox to\arrowwd{\hsize=0pt \vfil}%
}

\catcode`@=12

\def\Rzb#1{\hbox to 0pt{$#1$\hss}}
\def\Lzb#1{\hbox to 0pt{\hss$#1$}}
\def\Czb#1{\hbox to 0pt{\hss$#1$\hss}}
\def\rzb#1{\hbox to 0pt{$\scriptstyle#1$\hss}}
\def\lzb#1{\hbox to 0pt{\hss$\scriptstyle#1$}}
\def\czb#1{\hbox to 0pt{\hss$\scriptstyle#1$\hss}}
\def\ssRzb#1{\hbox to 0pt{$\scriptscriptstyle#1$\hss}}
\def\ssLzb#1{\hbox to 0pt{\hss$\scriptscriptstyle#1$}}
\def\ssCzb#1{\hbox to 0pt{\hss$\scriptscriptstyle#1$\hss}}

\def\Bzb#1#2{\hbox to0pt{\hbox to#2{\hss$#1$}\hss}}

\def\Bar#1{\setbox\n?box = \hbox{$#1$}
\kern .08\wd\n?box{\overline{\hbox to .84\wd\n?box{\vphantom {$#1$}\hss}}}
\kern .08\wd\n?box\kern -\wd\n?box\box\n?box}
\def\ctrlBar#1#2#3{\setbox\n?box=\hbox{$#1$}%
\n?dimen=\wd\n?box %
\advance\n?dimen by -#2pt\advance\n?dimen by -#3pt%
\kern #2pt\overline{\hbox to\n?dimen{\hss\vphantom{$#1$}}}%
\kern#3pt\kern-\wd\n?box\box\n?box}

\def\disjointunion{\hbox{$\perp\hskip -4pt\perp$}}


\def\Z{{\bf Z}}

\def\cy#1{\Z/{#1}\Z}

\newcount\hour
\newcount\minute
\newtoks\datetime
\newtoks\date
\def\?printdate{\datetime={\today}\date={\Today}}
\def\Today{\ifcase\month\or Jan.\or Feb.\or Mar.\or April\or May
\or June\or July\or Aug.\or Sept.\or Oct.\or Nov.\or Dec.\fi\space
\number\day, \number\year}
\def\today{\ifcase\month\or Jan.\or Feb.\or Mar.\or April\or May
\or June\or July\or Aug.\or Sept.\or Oct.\or Nov.\or Dec.\fi\space
\number\day, \number\year\ - %
{\count0=\time\divide\count0 by60 \hour=\count0%
\multiply\count0 by 60\minute=\time\advance\minute by-\count0%
\hour=\time\divide\hour by 60 \number\hour:\ifnum\minute<10 0\fi\number\minute}}
\?printdate%
\def\dateheaders{%
\headline{\hfil\ifnum\pageno>1 \sTOPr Printed \the\date\hfill\number\pageno.\fi\quad}%
\footline{\ifnum\pageno=1 \hbox to0pt{\sTOPr Printed \the\date\hss}\sTOPr\hfil 1.\fi\hfil}%
}%
\def\timeheaders{%
\headline{\hfil\ifnum\pageno>1 \sTOPr Printed \the\datetime\hfill\number\pageno.\fi\quad}%
\footline{\ifnum\pageno=1 \hbox to0pt{\sTOPr Printed \the\datetime\hss}\sTOPr\hfil 1.\fi\hfil}%
}%
\def\finalheaders{%
\headline{\hfil\ifnum\pageno>1 \hfill\number\pageno.\fi\quad}%
\footline{\ifnum\pageno=1 \hfil 1.\fi\hfil}%
}%
%

\newtoks\?ta\newtoks\?tb%
\long\def\append#1\to#2{\?ta{\\{#1}}\?tb=\expandafter{#2}%
\xdef#2{\the\?tb\the\?ta}}%
\def\lop#1\to#2{\expandafter\lopoff#1\lopoff#1#2}%
\long\def\lopoff\\#1#2\lopoff#3#4{\def#4{#1}\def#3{#2}}%
\def\cardinality#1\to#2{#2=0 \long\def\\##1{\advance#2 by1 }#1}%
\newbox\?title\def\?address{}%
\def\?thanks{}\def\?keywords{}\def\?subjclass{}%
\def\?translator{}%
\newbox\?abstract%
\def\?author{}\newcount\?authorcount%
\gdef\topmatter#1{#1}%
\font\sTOPl=cmsl8 %
\font\sTOPr=cmr8 %
\font\AbstractFont=cmr8 %
\font\TitleFont=cmbx12 at 20pt %
\font\AuthorFont=cmr12 %
\def\Titleskip{\vskip18pt}%
\def\Authorskip{\vskip10pt}%
\def\Topmatterskip{\vskip20pt}%
\gdef\endtopmatter{\centerline{\unhbox\?title}\Titleskip%
\cardinality\?author\to\?authorcount%
\def\\##1{##1\ifnum\?authorcount=2 \relax\ and \else%
\ifnum\?authorcount=1 \else , \fi\fi%
\global\advance\?authorcount by-1 }%
\centerline{\AuthorFont \?author}\Authorskip%
\itemitem{}{\unhbox\?abstract}\Topmatterskip%
\ifx\?subjclass\?empty\else%
\def\\{}\footnote{}{{\sTOPr1991 }%
{\sTOPl Mathematics Subject Classification}%
\ \sTOPr\?subjclass}\fi%
\ifx\?keywords\?empty\else%
\def\\{}\footnote{}{{\sTOPl Key words and phrases.}%
\ \sTOPr\?keywords}\fi%
\ifx\?thanks\?empty\else%
\def\\##1{##1\par}\footnote{}{\sTOPr\?thanks}\fi%
}%
\gdef\Title#1{\setbox\?title=\hbox{\TitleFont#1}}%
\gdef\AuthorInfo#1#2#3#4#5{%
\append#1\to\?author%
\def\?empty{}%
\def\y??y{\vskip4pt\noindent}%
\gdef\x??x{#2}\ifx\x??x\?empty\relax\else%
\edef\y??y{\y??y#2}\fi%
\gdef\x??x{#3}\ifx\x??x\?empty\relax\else%
\edef\y??y{\y??y\par{\noindent\sl email address:\ #3}}\fi%
\gdef\x??x{#4}\ifx\x??x\?empty\relax\else%
\edef\y??y{\y??y\par{\noindent\sl Current Address:\ #4}}\fi%
\gdef\x??x{#5}\ifx\x??x\?empty\relax\else%
\edef\y??y{\y??y\par{\noindent\sl URL:\ \tt#5}}\fi%
\expandafter\append\y??y\to\?address}%
\gdef\Thanks#1{\append#1\to\?thanks}%
\long\gdef\Abstract#1{\setbox\?abstract=\hbox{\AbstractFont Abstract:\ #1}}%
\global\let\?end=\end
\gdef\end{\ifx\?translator\?empty\relax\else%
\cardinality\?translator\to\?authorcount%
\def\\##1{##1\ifnum\?authorcount=2 \relax\ and \else%
\ifnum\?authorcount=1 \else , \fi\fi%
\global\advance\?authorcount by-1 }%
\par\noindent\hbox to\hsize{%
\hfill Translated by \?translator}%
\fi%
\long\def\\##1{\vskip4pt\noindent##1}\?address\?end}

\newcount\numberscheme
\newcount\label\label=0
\newif\ifexpandthmlabels\expandthmlabelsfalse
\newif\ifexpandbiblabels\expandbiblabelsfalse
\newif\ifexpandsymbollabels\expandsymbollabelsfalse
\newif\ifhyperlinked\hyperlinkedfalse
\numberscheme=1 

\def\check?#1#2{
\xdef\r?r{\expandafter\meaning\csname f?#1\endcsname}
\xdef\rr?r{\expandafter\meaning\csname relax\endcsname}%
\ifx\r?r\rr?r\relax
\xdef\r?r{\expandafter\meaning\csname f#1\endcsname}
{\ifx\r?r\rr?r\relax
XXX%
\else
\edef\x?x{\csname f#1\endcsname ?}\expandafter#2\x?x%
\fi}%
\else
\edef\x?x{ \csname f?#1\endcsname ?}\expandafter#2\x?x%
\fi%
}

\def\F?lakeI#1\ #2?{#1}
\def\F?lakeII#1\ #2?{#2}
\def\F?ull#1\ #2?{#1\ #2}

{{\catcode`#=11\xdef\Sharp{#}}}

\def\partII#1{\ifexpandthmlabels\relax\?addlabel{f#1}{}\fi%
\ifhyperlinked%
\check?{#1}{\F?lakeII}%
\else%
\check?{#1}{\F?lakeII}%
\fi}

\def\anchor#1{\ifexpandthmlabels\relax\?addlabel{f#1}{}\fi%
\ifhyperlinked%
{}\check?{#1}{\F?ull}%
\else%
\check?{#1}{\F?ull}%
\fi}

{{\catcode`~=11
\gdef\Tilde{~}
}}
\def\ExternalLink#1#2{\ifhyperlinked\relax
#2%
\else#2\fi}

\def\twothms#1#2{%
\edef\zm?{\csname f?#1\endcsname ?}%
\edef\xx?x{\expandafter\F?lakeI\zm?}%
\edef\zm?{\csname f?#2\endcsname ?}%
\edef\xy?x{\expandafter\F?lakeI\zm?}%
\ifx\xx?x\xy?x\relax\xy?x s \partII{#1} and \partII{#2}\else%
\csname f?#1\endcsname\ and \csname f?#2\endcsname\fi}

\def\stufflabel{\ifnum\numberscheme=1 \number\sectionnumber.\number\label
\else\number\label\fi}

\def\tmplabel{xx}


\def\autonewlabel#1#2#3{\global\advance\label by 1 %
\??autonewlabel{#1}{#2}{#3}}

\def\rfc#1#2?{#2}
\def\macroname#1{%
\expandafter\rfc\string#1 ?}

\def\??autonewlabel#1#2#3{%
\def\xxcc{#1}\ifx\xxcc\tmplabel\relax\else
\xdef\r?r{\expandafter\meaning\csname f?#1\endcsname}
\xdef\rr?r{\expandafter\meaning\csname relax\endcsname}%
\ifx\r?r\rr?r\relax%
\else%
\xdef\r?r{\csname f?#1\endcsname}
\let\nl?c=\newlinechar%
\newlinechar=`\^^J%
\errmessage{\string f#1\space is already defined as^^J\r?r}%
\let\newlinechar=\nl?c%
\fi\fi%
\def\r?r{#3}%
\ifx\r?r\?empty\relax\edef\r?r{\stufflabel}\fi
\ifexpandthmlabels\relax
\expandafter\xdef\csname f#1\endcsname{\ifhyperlinked%
\?addlabel{f#1}{#2\ \r?r}%
\else%
\?addlabel{f#1}{#2\ \r?r}\fi%
}%
\expandafter\xdef\csname nlf#1\endcsname{\?addlabel{f#1}{#2\ \r?r}}%
\else%
\expandafter\xdef\csname f#1\endcsname{\ifhyperlinked%
#2\ \r?r%
\else%
#2\ \r?r\fi}%
\expandafter\xdef\csname nlf#1\endcsname{#2\ \r?r}%
\fi%
\expandafter\xdef\csname f?#1\endcsname{#2\ \r?r}
\ifwritelabellist%
\immediate\write\labellist{\def\string\f#1{#2\ \r?r}}
\fi
}

\def\?addlabel#1#2{\setbox0=\hbox{\raise 8pt\hbox to 0pt{%
$\scriptscriptstyle\backslash${\smallfont #1}\hss}{#2}}\unhbox0}

\def\newT#1{\autonewlabel{#1}{\?thm?}{}}

\def\newL#1{\autonewlabel{#1}{\?lem?}{}}

\def\newP#1{\autonewlabel{#1}{\?pro?}{}}

\def\newR#1{\autonewlabel{#1}{\?rem?}{}}



\def\?thm?{Theorem}
\def\?cor?{Corollary}
\def\?lem?{Lemma}
\def\?def?{Definition}
\def\?pro?{Proposition}
\def\?for?{Formula}
\def\?exa?{Example}
\def\?dia?{Diagram}
\def\?rem?{Remark}
\def\?res?{Remarks}
\def\ab??X#1\ #2?{\def\x??x{#1}%
\ifx\x??x\?thm? Thm.\else%
\ifx\x??x\?cor? Cor.\else
\ifx\x??x\?lem? Lem.\else%
\ifx\x??x\?def? Def.\else%
\ifx\x??x\?pro? Prop.\else%
\ifx\x??x\?for? Form.\else%
\ifx\x??x\?exa? Ex.\else%
\ifx\x??x\?dia? Diag.\else%
\ifx\x??x\?rem? Rem.\else%
\ifx\x??x\?res? Rem.'s\else%
\x??x%
\fi\fi\fi\fi\fi\fi\fi\fi\fi\fi\ #2}

\def\ab?X#1?{\edef\x??x{#1}\show\x??x\expandafter\ab??X\x??x}

\newcount\referencecount\referencecount=0
\newif\ifusereferencecount\usereferencecountfalse


\def\ref?list{}
\def\newref#1#2#3{\global\advance\referencecount by 1 %
\xdef\r?r{\expandafter\meaning\csname b?#1\endcsname}
\xdef\rr?r{\expandafter\meaning\csname relax\endcsname}%
\ifx\r?r\rr?r\relax\else%
\xdef\r?r{\csname b?#1\endcsname}%
\let\nl?c=\newlinechar%
\newlinechar=`\^^J%
\errmessage{\string b#1\space is already defined as^^J\r?r}%
\let\newlinechar=\nl?c%
\fi%
\ifusereferencecount\relax\def\r?r{\number\referencecount}\else%
\def\r?r{#2}\ifx\r?r\?empty\relax\def\r?r{\number\referencecount}%
\fi\fi%
\ifwritelabellist\immediate\write\labellist{%
\def\string\b #1{\r?r}}\fi%
\ifexpandbiblabels\relax%
\expandafter\xdef\csname b#1\endcsname{\ifhyperlinked%
\?addlabel{b#1}{\r?r}%
\else%
\?addlabel{b#1}{\r?r}\fi}%
\else%
\expandafter\xdef\csname b#1\endcsname{\ifhyperlinked%
\r?r\else\r?r\fi}
\fi%
\ifexpandbiblabels\relax%
\expandafter\xdef\csname Ab#1\endcsname{\ifhyperlinked%
\?addlabel{b#1}{\r?r}%
\else%
\?addlabel{b#1}{\r?r}\fi}%
\else%
\expandafter\xdef\csname Ab#1\endcsname{\ifhyperlinked%
\r?r\else\r?r\fi}
\fi%
\expandafter\xdef\csname b?#1\endcsname{\r?r}%
\expandafter\append\csname Ab#1\endcsname\to\ref?list%
}

\def\cite[#1]{[#1]}
\def\Cite[#1,#2]{[#1,#2]}

\newwrite\labellist
\newread\tmpinput

\newif\ifwritelabellist\global\writelabellistfalse
\font\smallfont=cmr6

\def\usedotlabelfile{\input \jobname.label
}

\long\def\?loop#1\?repeat{\def\?body{#1}\?iterate}
\long\def\?iterate{\?body\let\?next=\?iterate\fi\?next}

\def\rewritedotlabelfile{%
\ifwritelabellist%
\immediate\closeout\labellist%
\immediate\openout\labellist=\jobname.label%
\immediate\openin\tmpinput=tmp.label%
\catcode`\\=12 %
\def\r?r{\par}%
\?loop%
\ifeof\tmpinput\let\?next=\relax\else
\read\tmpinput to\l??ne
\ifx\l??ne\r?r\else
\immediate\write\labellist{\l??ne}\fi
\?repeat
\catcode`\\=0 %
\fi}

\catcode`?=12 \relax%
\numberscheme=0
\newref{BS}{}{}
\newref{CL}{}{}
\newref{GH}{}{}
\newref{R}{}{}
\newref{TT}{}{}

\def\ucc#1{r_{conj}(#1)}


\topmatter{\finalheaders\font\TitleFont=cmbx12}
\Title{Codimension one spheres which 
are null homotopic.}
\AuthorInfo{Laurence R.~Taylor\Thanks{Partially supported by the N.S.F.}}
{Department of Mathematics, University of Notre Dame, Notre Dame, IN 46556}
{taylor.2@nd.edu}{}{\ExternalLink{http://www.nd.edu/\Tilde taylor}
{http://www.nd.edu/\Tilde taylor}}
\endtopmatter

Grove and Halperin \cite[\bGH] introduced a notion of {\sl
taut immersions}.
Terng and Thorbergsson \cite[\bTT] give a slightly different definition and
showed that taut immersions are
a simultaneous generalization of taut immersions of manifolds into
Euclidean spaces or spheres, and some interesting embeddings constructed by Bott
and Samelson \cite[\bBS].
They go on to prove many theorems about such immersions.
One particularly intriguing result, Theorem 6.25, concerned
codimension one, null homotopic, tautly embedded
spheres.
Using a result of Ruberman, \cite[\bR], they proved
in many cases that this sphere had to be a distance
sphere, that is, the image of a standard sphere 
under the exponential map, a generalization of a theorem of
Chern and Lashof \cite[\bCL].
Here we observe that the methods of Terng{--}Thorbergsson and
Ruberman suffice to classify tautly immersed, null homotopic,
codimension one spheres.
Informally one may say that the examples produced by Terng and Thorbergsson
in \cite[\bTT] are all that there are.
Precisely, we have

\medskip
\newT{TT}
\Thm{\anchor{TT}}{
Let $N^n$ be a complete Riemannian manifold and let
$\phi\colon S\to N$ be a tautly immersed, codimension one
sphere which is null homotopic.
Then there exists a taut point $q\in N$ and a real number $l$, $0<l<\ucc{q}$, 
such that $S$ is the image of the sphere of radius $l$ in $T_qN$
under the exponential map.
Here $\ucc{q}$ denotes the conjugate radius at 
the point $q\in N$.}

\Def{Remark}{
If we fix $q\in N$ and a real number $l$, $0<l<\ucc{q}$, then
the image of the sphere of radius $l$ under the exponential
map from $T_qN$ is a taut immersion if and only if $q$ is taut.
This is Theorem 6.23, p.~219, of Terng and Thorbergsson \cite[\bTT].
It follows from the definition of conjugate radius that $exp_q\colon T_qN\to N$
is an immersion when restricted to the ball of radius $l$ and from Theorem 6.23
\cite[\bTT] it follows that the spheres
parallel to the boundary are all taut.}

To explain the result of Ruberman eluded to above,
we introduce some notation.
A locally flat, codimension one sphere
in $M$ which disconnects $M$ displays $M$ as 
a connected sum, $M=M_1\# M_2$.
We call this sphere a {\sl connected sum sphere}.
If we assume a codimension one, locally flat
sphere is null homotopic, it
must be a connected sum sphere.
One obvious way to get null homotopic, connected sum spheres 
is for at least one of $M_1$ or $M_2$ to be a homotopy sphere. 
It is natural to wonder if there are any others.

Ruberman \cite[\bR] gives examples of null homotopic 
connected sum spheres for which neither $M_1$
nor $M_2$ is a homotopy sphere.
He also gives a nearly necessary and sufficient criterion 
to have
a null homotopic, connected sum sphere, \Cite[\bR,Thm.~1].
We observe that a
more vigorous use of Ruberman's main lemma provides a
complete characterization.
This characterization shows that Ruberman's list of
examples is complete.
The characterization is given by the following theorem.

\newT{R}
\Thm{\anchor{R}}{
A connected sum sphere 
in $M$, a connected, paracompact, Hausdorff, $m${--}dimensional  
topological manifold, with or without boundary, 
is null homotopic 
if and only if either (1) or (2) below holds, where
$M=M_1\#M_2$.
\Itemitem{(1)} At least one of $M_1$ or $M_2$ 
is a homotopy $m${--}sphere.
\Itemitem{(2)} Neither $M_1$ nor $M_2$ is a homotopy sphere,
but there exist relatively prime integers $k_1$ and $k_2$,
both greater than $1$, so that each $M_i$ is a
simply connected
$\Z[\Frac{1}{k_i}]${--}homology $m${--}sphere.
}

\Def{Remarks}{
We repeat Ruberman's remark that case (2) can not hold
if $m\leq 4$, but there are infinitely many examples of (2) 
in each dimension $\geq 5$.
We emphasize that case (1) must occur unless
$M$ is a compact, simply connected manifold without
boundary which is a rational homology $m${--}sphere with
torsion in its homology containing at least two distinct primes.
We also remark that Ruberman already proved \fR\ in the case $M$ is simply connected
and this is the only case we use in the proof of \fTT.}

The author would like to thank D.~Ruberman and 
G.~Thorbergsson for useful conversations.

\section{{The proof of \nlfR.}}

To state Ruberman's main lemma requires further notation
and conventions.
To form the connected sum requires embeddings
$\iota_i\colon B^m\to M_i$.
Then $M_1\# M_2$ is formed by removing the interiors of
$\iota_i(B^m)$ and identifying the boundaries.
If $M_1\vee M_2$ denotes the wedge $M_1$ and $M_2$
then there is a natural map
$M_1\# M_2\to M_1\vee M_2$ obtained by pinching 
$i(S^{m-1})$ to a point and identifying $M_i/\iota_i(B^m)$
with $M_i$ via an inverse to the map
$M_i\to M_i/\iota_i(B^m)$ which is the identity outside a 
small neighborhood of $B^m$.
Let $\hat B_i\subset M_i$ be a larger ball containing
$\iota_i(B^m)$ so that 
$\overline{\hat B_i-\iota_i(B_i)}=S^{m-1}\times[0,1]$.
There are maps $M_i\to S^m$ which pinch all of 
$M_i-\hat B_i$ to a point.
Composition yields a map
$M_1\# M_2\to M_1\vee M_2\to S^m\vee S^m$.
Orient $S^{m-1}=\partial B^m$ and use the inward normal 
last of $\partial B^m$ to orient
the balls in $M_1$ and $M_2$.
This orients the two $S^m$ as well.

Now suppose the connected sum sphere is null homotopic.
This means that $i$ extends to a map of the disk
$I\colon D^m\to M_1\# M_2$.
The map $I$ yields two maps $I_i\colon S^m\to M_i$.
Define the {\sl degree\/} of the map $I_i$ as the
degree of the composite $S^m\RA{I_i}M_i\to S^m$
where we orient the first $S^m$ from $D^m$ using the
inward normal last from the orientation on $S^{m-1}$.
Ruberman's main lemma then becomes

\newL{RML}
\Thm{\anchor{RML}. {\rm (Ruberman \Cite[\bR,Lemma~3])}}{
With notation and conventions as above, 
$$\deg I_1+\deg I_2=1\ .$$}

\noindent
For completeness we give a proof.
The following diagram commutes
$$\Matrix{M_1\# M_2&\to&M_1\vee M_2&\to&M_i\CR
\big\downarrow&&\big\downarrow&&\big\downarrow\CR
S^m\# S^m&\to&S^m\vee S^m&\to&S^m\CR
}$$
It follows that it suffices to prove the result for
the case $M_1=M_2=S^m$.
Orient $S^m\# S^m$ so that the map
$\Z\simeq H_m(S^m\#S^m)\to H_m(S^m\vee S^m)\simeq\Z\oplus \Z$
sends $1$ to $(1,-1)$.
Given two null homotopies, $I, J\colon D^m\to S^m\# S^m$
of the connected sum sphere,
they differ by an element in $\pi_m(S^m\# S^m)$ and the
map $\Z\simeq\pi_m(S^m\# S^m)\to H_m(S^m\vee S^m;\Z)
\simeq\Z\oplus\Z$
sends $1$ to $(1,-1)$.
Hence the sum of the degrees
is constant.
The evident null homotopy which collapses the connected sum
sphere in $S^m\# S^m$
using the first $S^m$ has degrees $(1,0)$, so
that difference is $1$.\Qed

Two further lemmas will be required.

\newL{RSC}
\Thm{\anchor{RSC}}{If $X$ is a simply connected,
$\Z[\Frac{1}{k}]${--}acyclic space, then $\pi_{m-1}(X)$
is $k${--}torsion, $m\geq 2$.}

As Ruberman observes, this is an elementary consequence
of Serre's mod ${\cal C}$ theory.
He also gives a proof of this next lemma as his Lemma 2.

\newL{RDOM}
\Thm{\anchor{RDOM}}{
If $f\colon S^m\to M^m$, $m>1$, has degree $k$, $k\neq 0$,then
$\pi_1(M)$ is finite of order dividing $k$.
It further follows that $M$ is a
$\Z[\Frac{1}{k}]${--}homology $m${--}sphere.}

\newR{IR}
\Def{\anchor{IR}}{
Note that a simply connected homology sphere is a 
homotopy sphere.
If $k=\pm1$ in \fRDOM, $M$ is a homotopy sphere 
even if $m=1$.
If $M$ is not closed, compact 
and orientable, then automatically $k=0$.}

\medskip
First we prove that the connected sum sphere in case (1) or in
case (2) is null homotopic.
The proof presented is essentially Ruberman's.
In case (1), the connected sum sphere is null homotopic
because $M_i-D^m$ is a homotopy disk if $M_i$ is a
homotopy sphere.
In case (2), recall that $m\geq 5$.
Note that $M_i-D^m$ is $\Z[\Frac{1}{k_i}]${--}acyclic
and the map $S^{m-1}\RA{\ i\ }\ M_1\# M_2$ factors through 
each $M_i-D^m$.
By \fRSC, the class in $\pi_{m-1}(M_1\# M_2)$ represented
by the connected sum sphere is both $k_1$ torsion 
and $k_2$ torsion.
Since $k_1$ and $k_2$ are relatively prime, 
$i$ is null homotopic.

{\sl Now suppose the connected sum sphere is null homotopic 
and that neither side is a homotopy sphere.
We must show case (2) holds.\/}

\medskip\noindent{\bf Step 1.}
{\sl We have $m>1$ and there exist relatively prime numbers
$k_1$ and $k_2$, both greater than $1$ such that
each $M_i$ is a $Z[\Frac{1}{k_i}]${--}homology
$m${--}sphere.
The fundamental group of $M_i$ is finite of
order dividing $k_i$.}

\pf
From \fRDOM, \fIR\ and
our assumption that neither side is a homotopy
sphere, it follows that neither $\deg I_i$ can be $\pm1$.
It follows from \fRML\ that neither $\deg I_i$ can be
$0$ either.
It follows from \fIR\ that $M_i$ is closed, compact, 
orientable.
If $m=1$, $M_i=S^1$ which is a contradiction.
Hence $m>1$.
Let $k_i$ be the absolute value of $\deg I_i$ and note that
each $k_i$ is greater than $1$.
It follows from \fRML\ that 
the $k_i$ are relatively prime.
It follows from \fRDOM\ that each 
$M_i$ has finite fundamental
group of order dividing $k_i$ 
and is a $\Z[\Frac{1}{k_i}]${--}homology 
$m${--}sphere.\QED

\smallskip If $M$ is a rational homology sphere with finite
fundamental group, let $r(M)$ denote the order of $\pi_1(M)$ 
and let $\ell(M)$ denote the order of the direct sum of the
torsion subgroups of the homology.
Step 1 can be summarized as

\Thm{Summary}{
If the connected sum sphere in $M^m=M_1\# M_2$ is null homotopic
with $M$ connected and neither $M_i$ a homotopy sphere,
then each $M_i$ is a rational homology $m${--}sphere,
$\bigl(\ell(M_1),\ell(M_2)\bigr)=1$, and
$\bigl(r(M_1),r(M_2)\bigr)=1$.
Since neither $M_i$ is a homotopy sphere, 
$m>1$ and
$\ell(M_i)\cdot r(M_i)>1$.}

Consider the following construction.
Let $\widetilde M_i$ denote the universal cover of $M_i$
and let $Y_i$ be the connected sum of
$\widetilde M_i$ and $(r(M_i)-1)$ copies of $M_{3-i}$.

\medskip\noindent{\bf Step 2.} {\sl $Y_i$ is not a homotopy
sphere.}

\pf
Since $M_{3-i}$ is not a homotopy sphere, neither is
$Y_i$ if $r(M_i)>1$.
But if $r(M_i)=1$, $Y_i=M_i$ is not a homotopy
sphere either.\QED

The next two formulae follow from the Mayer{--}Vietoris
and the van{--}Kampen theorems.
\NumAlign{
\ell(Y_i) =& \ell(\widetilde M_i)\cdot 
\ell(M_{3-i})^{r(M_i)-1}\CR{(5)}
\pi_1(Y_i)=&\pi_1(M_{3-i})\ast\cdots\ast\pi_1(M_{3-i})
\qquad r(M_i)-1\ {\rm times.}\CR{(6)}}

\medskip\noindent{\bf Step 3.} {\sl The connected sum
sphere in $Y_i\# M_{3-i}$ is null homotopic.}

\pf
The manifold $Y_i\# M_{3-i}$ is the total space of a
cover of $M_i\#M_{3-i}$:
$\pi\colon Y_i\#M_{3-i}\to M_i\# M_{3-i}$
denotes the covering map.
Check that $\pi$ of the connected sum sphere in 
$Y_i\# M_{3-i}$ is the connected sum sphere in 
$M_i\# M_{3-i}$.
Since $m\geq2$, the induced map 
$\pi_\ast\colon\pi_{m-1}(Y_i\# M_{3-i})
\to \pi_{m-1}(M_i\# M_{3-i})$
is an injection.\QED

\medskip\noindent{\bf The Final Step:} {\sl 
$\pi_1(M_i)=0$: i.e.\ $r(M_i)=1$.}

\pf 
Suppose $r(M_i)>1$.
By Steps 2 and 3, the Summary applies to
$Y_i\# M_{3-i}$.
Hence $\ell(Y_i)$ and $\ell(M_{3-i})$ are relatively prime,
so $\ell(M_{3-i})=1$ by (5).
It follows from the Summary then that $r(M_{3-i})>1$.
Since $\pi_1(Y_i)$ is finite, (6) implies $r(M_i)=2$ 
and hence $r(Y_i)=r(M_{3-i})$. 
This is a contradiction since $\bigl(r(Y_i),r(M_{3-i})\bigr)=1$.
Hence $r(M_i)=1$.\QED

\section{The proof of \nlfTT.}

We begin with three results we need in the proof.
First observe that more mileage is available from the
proof of Theorem 2.5, p.~188 of \cite[\bTT].

\newT{emb}
\Thm{\anchor{emb}}{
Let $\phi\colon M\to N$ be a taut immersion and fix a point
$x\in M$.
Assume $M$ is connected.
Let $\pi\colon \widetilde N\to N$ be the cover for which
there is a choice of $y\in\widetilde N$ such that
the subgroup
$\pi_\ast\bigl(\pi_1(\widetilde N, y)\bigr)\subset 
\pi_1(N,\phi(x))$ equals the subgroup
$\phi_\ast\bigl(\pi_1(M,x)\bigr)\subset 
\pi_1(N,\phi(x))$.
Let $\tilde\phi\colon M\to\widetilde N$ be the unique lift
with $\tilde\phi(x)=y$.
Then $\tilde\phi$ is an embedding.}

\pf
Note that $P(\widetilde N,\tilde\phi\times q)$ is connected since
$M$ is and $\tilde\phi$ induces an isomorphism on $\pi_1$.
Next note that $P(\widetilde N,\tilde\phi\times q)$ is a 
component of $P\bigl(N,\phi\times \pi(q)\bigr)$.
If we choose $q\in\widetilde N$ as Terng and Thorbergsson do 
on page 188, then the energy function is a perfect Morse function on 
$P(\widetilde N,\tilde\phi\times q)$.
The result follows just as in \cite[\bTT].\QED

\newP{cL}
\Thm{\anchor{cL}}{
Let $\phi\colon M\to N$ be a map, let $\pi\colon \widetilde N\to N$
be a cover, and let $\tilde\phi\colon M\to \widetilde N$
be a map covering $\phi$.
Then $\phi$ is a taut immersion if and only if $\tilde\phi$ is.}

\pf
It follows from the structure of covering spaces that $\phi$ is an
immersion if and only if $\tilde\phi$ is.
Next note that an immersion $\psi\colon X\to Y$ is taut if and
only if the energy function on the path space $P(Y,\psi\times q)$
is a perfect Morse{--}Bott function for all points $q\in Y$.
One direction is Theorem 2.9, p.~193, in \cite[\bTT].
The other direction follows from the remark on page 183 of \cite[\bTT]
that the energy function on the path space is Morse at all non{--}focal
points, plus the observation that a proper Morse{--}Bott function which
is Morse is a proper Morse function.

Now the path space for $\phi$
is the disjoint union of various path spaces
$$P(N,\phi\times p)=\mathop{\disjointunion}_{p_\alpha\in \pi^{-1}(p)}
P(\widetilde N,\tilde\phi\times p_\alpha)\ .$$
Since the energy function on each $P(\widetilde N,\tilde\phi\times p_\alpha)$
is just the restriction of the energy function on
$P(N,\phi\times p)$, the result is immediate.
\QED

\Def{Remark}{G.~Thorbergsson has shown the author an elementary proof
of \fcL\ that does not need Theorem 2.9.}

Finally we need a generalization of Theorem 6.25 \cite[\bTT].

\newP{mtb}
\Thm{\anchor{mtb}}{
Let $\psi\colon S\to M$ be a taut embedding and assume
$\psi(S)$ bounds a simply connected, $\cy2$ homology ball $B\subset M$.
Then $B$ is a ball and there exists a taut point $p\in B\subset M$ and a real
number $l>0$ such that $\psi(S)$ is the image of the sphere of radius $l$ in
$T_pM$ and $l$ is less than $\ucc{p}$.
}

\pf
We indicate the changes needed from the proof of Theorem 6.25 on pages
219 and 220 of \cite[\bTT].
Following Terng and Thorbergsson, pick $p\in B\subset M$ so that
there is a geodesic perpendicular to $\psi(S)$ with first 
focal point $p$ at a distance $l$ from $\psi(S)$.
We further assume $l$ is minimal with this property. 
Observe that the space $P(B,\psi\times p)$
fits into a fibration sequence
$$\Omega B\to P(B,\psi\times p)\ \RA{\ \pi\ }\ S\ .$$
Since $B$ is simply connected and mod $2$ acyclic,
it follows from the Serre spectral sequence that
$\Omega B$ is also mod $2$ acyclic.
It then follows from the Serre spectral sequence that
$\pi_\ast$ induces an isomorphism in mod $2$ homology.
Hence $P(B,\psi\times p)$ has no homology below dimension $n-1$.
Argue just as in Terng and Thorbergsson that the multiplicity of the focal point
$p$ is $n-1$, so all the geodesics starting perpendicularly to $\psi(S)$
meet at a distance $l$ in the point $p$.
Furthermore, all of these geodesics minimize
distance between $p$ and $\psi(S)$.
Hence the conjugate radius at $p$ is bigger then $l$.
Furthermore no two of these geodesics intersect in $B$
since they minimize distance, so
$B$ is the ball of radius $l$ centered at $p$.
It follows from \Cite[\bTT,Thm.~6.23, p.~219] that $p$ is taut.\QED

To fix notation for the proof of \fTT, let $\pi\colon\widetilde N\to N$
be the universal cover.
Since $\phi\colon S^{n-1}=S\to N$ is null homotopic, $\phi$ lifts
to the universal cover via $\tilde\phi\colon S\to \widetilde N$.
By \femb, $\tilde\phi$ is an embedding and it remains null homotopic.
Hence by \fR, there is a simply connected, $\cy2$ homology ball, 
$B\subset\widetilde N$, with boundary $\partial B=\tilde\phi(S)$.
The embedding $\tilde\phi$ is taut by \fcL.
Applying \fmtb, we now have a taut point 
$\tilde q\in\widetilde N$ and a real number $l$ such that 
$\tilde\phi(S)$ is the distance sphere of radius $l$ centered
at $\tilde q$ with $l<\ucc{\tilde q}=\ucc{q}$, where $q$
denotes the image of $\tilde q$ in $N$.
It follows that
$\phi(S)$ is the image of the sphere of radius $l$ in $T_qN$
under the exponential map.
Finally apply \fcL\ again: since $\tilde q\in\widetilde N$
is taut, $q$ is taut
in $N$ as required by \fTT.

\medskip
\Def{Remark}{The hypothesis in \fTT\ that the immersion $\phi$ is null
homotopic is sometimes unnecessary.
By \Cite[\bTT,Thm.~2.5], if $\pi_1(N)=0$ then $\phi$ is an
embedding and in many cases, the $\cy2$ cohomology ring will force
one of the pieces in the connected sum decomposition to be a simply connected 
$\cy2$ homology ball.
Then \fmtb\ will show that $\phi$ is null homotopic.
Examples include all simply connected manifolds with the $\cy2$ cohomology of 
the rank $1$ symmetric spaces, $H^\ast(S^n;\cy2)$, $H^\ast(CP^n;\cy2)$,
$H^\ast(HP^n;\cy2)$ and the $\cy2$ cohomology of the Cayley projective plane.
Products of such spaces are also examples.}

\Def{Remark}{One can generalize Ruberman's theorem to show that if the
order of $\phi\in\pi_{n-1}(N)$ is odd, then the embedding $\tilde\phi$
bounds a simply connected $\cy2$ homology ball, so again \fmtb\ will
show that $\phi$ is null homotopic.}

\medskip
\penalty-1000
\noindent{\bf Bibliography}
\newdimen\labellength
\newdimen\bXw
\newdimen\bXXW
\bXXW=4pt

\def\bibitem#1{\vskip 6pt%
\vfil\penalty-100\vfilneg%
\hangindent=\labellength%
\hangafter=1 %
\setbox0=\hbox {[\ {\bf#1}\ ]}%
\bXw=\labellength\advance\bXw by-\wd0\advance\bXw by-\bXXW
\noindent\hskip\bXw [\ {\bf#1}\ ]\hskip\bXXW}%
\setbox0=\hbox{\bf [99\ ]}
\labellength=\wd0
\global\advance\labellength by 4 true pt
\penalty10000
\bibitem{\AbBS}%
R.~Bott and H.~Samelson, {\sl Applications of the theory of {M}orse to
  symmetric spaces}, Amer. J. Math. {\bf 80} (1958), 964--1029, Corrections in
  {\bf 83} (1961) 207-208.

\bibitem{\AbCL}%
S.-S. Chern and R.~K. Lashof, {\sl On the total curvature of immersed
  manifolds}, Amer. J. Math. {\bf 79} (1957), 306--318.

\bibitem{\AbGH}%
K.~Grove and S.~Halperin, {\sl Elliptic isometries, condition $({\rm
  {C}})$ and proper maps}, Arch. Math. (Basel) {\bf 56} (1991), 288--299.
\bibitem{\AbR}%
D.~Ruberman, {\sl Null homotopic spheres of codimension one}, ``Tight and Taut
  Submanifolds'' (T.~E. Cecil and S.-S. Chern, eds.), MSRI Publications,
  vol.~32, Cambridge Univ. Press, Cambridge, 1997, pp.~229--232.

\bibitem{\AbTT}%
C.-L. Terng and G.~Thorbergsson, {\sl Taut immersions into complete Riemannian
  manifolds}, ``Tight and Taut Submanifolds'' (T.~E. Cecil and S.-S. Chern,
  eds.), MSRI Publications, vol.~32, Cambridge Univ. Press, Cambridge, 1997,
  pp.~181--228.

\end